\documentclass[titlepage,12pt]{article}
\usepackage{amssymb}
\usepackage{amsfonts}
\usepackage{amsmath}
\begin{document}

{\bf \Large Gap in Nonlinear Equivalence for \\ \\ Numerical Methods for PDEs} \\ \\

{\bf Elemer E Rosinger} \\
Department of Mathematics \\
and Applied Mathematics \\
University of Pretoria \\
Pretoria \\
0002 South Africa \\
eerosinger@hotmail.com \\ \\

{\bf Abstract} \\

For a large class of nonlinear evolution PDEs, and more generally, of nonlinear semigroups, as well as their approximating numerical methods, two rather natural stability type convergence conditions are given, one being necessary, while the other is sufficient. The gap between these two stability conditions is analyzed, thus leading to a general nonlinear equivalence between stability and convergence. \\  \\

{\bf 1. The General Setup} \\

The study of linear and nonlinear evolution systems of PDEs, with possibly associated initial and/or boundary problems can, as is well known, be dealt with in the more general framework of semigroups depending on a continuous parameter which represents time. Here, the study of the nonlinear equivalence for numerical methods approximating exact solutions of nonlinear evolution PDEs will be dealt with in this more general framework, namely, of numerical
methods approximating rather general nonlinear semigroups. \\

{\bf Definition 1.1.} \\

Given a normed vector space $( X, ||~|| )$. By a {\it nonlinear semigroup} on $X$ we mean any
family of mappings \\

(1.1)~~~ $ E ( t ) : X_t \longrightarrow X,~~~ t \in [ 0, \infty ) $ \\

where \\

(1.2)~~~ $ X_{t + s} \subseteq X_t \subseteq X,~~~ t, s \in [ 0, \infty ) $ \\

(1.3)~~~ $ E ( 0 ) = id_{X_0} $ \\

while for every $t, s \in [ 0, \infty )$ we have the commutative diagram \\

\begin{math}
\setlength{\unitlength}{0.2cm}
\thicklines
\begin{picture}(60,21)

\put(10,17){$X_{t + s}$}
\put(27,19){$E ( t + s )$}
\put(15,17.5){\vector(1,0){31}}
\put(47,17){$X$}
\put(0,7){$(1.4)$}
\put(15,15){\vector(1,-1){13}}
\put(16,7){$E ( t )$}
\put(29,0){$X_s$}
\put(32,2){\vector(1,1){13}}
\put(40,7){$E ( s )$}

\end{picture}
\end{math}

\hfill $\Box$ \\

An immediate consequence of (1.4) is that (1.1) takes the stronger form \\

(1.5)~~~ $  E ( t ) : X_t \longrightarrow X_0,~~~ t \in [ 0, \infty ) $ \\

The nonlinear semigroup (1.1) - (1.5) is called {\it well posed}, if and only if the mapping \\

(1.6)~~~ $ Z \ni ( t, u ) \longmapsto E ( t ) u \in X $ \\

is continuous, where \\

(1.7)~~~ $ Z = \bigcup_{t \in [ 0, \infty )}~ \{ t \} \times X_ t $ \\ 

{\bf Remark 1.1.} \\

The reason for {\it not} asking that $X_t = X$, for $t \in [ 0, \infty )$, is to allow for "finite blow up time" which is a frequent phenomenon in the nonlinear case, as illustrated by a simple example like the semigroup $E ( t )$ of solutions of the ODE \\

$~~~ U\,' ( t ) = U^2 ( t ),~ U ( 0 ) = u,~ t \in [ 0, \infty ) $ \\

in which case we can take $( X, ||~|| ) = \mathbb{R}$ and $X_t = ( - \infty, 1 / t )$, for $t \in ( 0, \infty )$, while $X_0 = X$, since for $u \in X_t$ we have $U ( t ) = u / ( 1 - u t )$. In this case $E ( t ) u = U ( t )$, for $t \in [ 0, \infty )$. \\ \\

{\bf 2. Three Versions of Stability} \\

The reason stability conditions for numerical methods solving PDEs are so important from practical point of view is very simple : such conditions are formulated exclusively in terms of the respective numerical methods, thus their verification does not in any way involve any knowledge of the exact solutions of the respective PDEs, and instead, it can be performed based alone on the given numerical methods. \\

The classical example in this regard is the Lax-Richtmyer Linear Equivalence Theorem, [Lax-Richtmyer], seen as the fundamental theorem of linear numerical analysis, and which says that the stability of a numerical method is both necessary and sufficient for its convergence to the exact solution of the respective linear evolution PDE, provided that the PDE is well posed, and the numerical method is consistent with that PDE. \\
Here of course, the well posedness of the PDE does not depend on any numerical method, being but an intrinsic property of the PDE under consideration. As for the consistency of the numerical method with the given PDE, that can usually be established based on suitable smoothness assumptions on the exact solutions of the PDE, without any more specific knowledge of such exact solutions. \\
And then, once well posedness and consistency have been established, the crucial issue of convergence of the numerical method to the unknown exact solution of the given PDE is simply identical with the issue of the stability of that numerical method, stability which, as mentioned, is an intrinsic property of the numerical method, thus it can be established without any further knowledge about the respective PDE or about its exact solutions. \\

One of the main {\it troubles} with the Lax-Richtmyer Linear Equivalence Theorem, Rosinger [1-9], is that its proof depends essentially on the Principle of Uniform Boundedness of Linear Operators in Banach Spaces, a principle which simply cannot have a general enough nonlinear correspondent. Therefore, general enough {\it nonlinear} extensions of the Lax-Richtmyer Linear Equivalence Theorem are hard to come by. \\
Here however, it is important to note that in the proof of the Lax-Richtmyer Linear Equivalence Theorem, the implication "stable $\Longrightarrow$ convergent" is trivial, and the hard part is the proof of the converse implication "convergent $\Longrightarrow$ stable", where the mentioned essentially linear principle is made use of. \\
Therefore, one may expect that in any proper and wide enough nonlinear generalization of the Lax-Richtmyer Linear Equivalence Theorem the difficult part would be the implication "convergent $\Longrightarrow$ stable". \\ 
In other words, finding {\it necessary} conditions of nonlinear stability may appear to be the main challenge. \\

Here, following Rosinger [2], we present {\it three} general nonlinear stability conditions which, as seen in the sequel, are closely related. \\

{\bf Definition 2.1.} \\

Given a well posed nonlinear semigroup (1.1) - (1.7), the numerical methods associated with it will be {\it continuous} mappings \\

(2.1)~~~ $ [ 0, \infty ) \times X \ni ( \Delta t, u ) \longmapsto C_{\Delta t} u \in X $ \\

which satisfy \\

(2.2)~~~ $ C_0 = id_X $

\hfill $\Box$ \\

In the sequel, we need the following auxiliary notion related to (1.7). A subset $Z\,' \subseteq Z$ is called {\it regular}, if and only if the following
three conditions are satisfied \\

(2.3)~~~ $ Z\,' ~~\mbox{is closed} $ \\

(2.4)~~~ $ E ( t ) X\,'_{t + s} \subseteq X\,'_s,~~~ t, s \in [ 0, \infty ) $ \\

(2.5)~~~ $ C_{\Delta t} X\,'_{t + \Delta t} \subseteq X\,'_t,~~~ t, \Delta t \in [ 0, \infty ) $ \\

where we denoted \\

(2.6)~~~ $ X\,'_t = \{\, u \in X_t ~|~ ( t, u ) \in Z\,' \,\},~~~ t \in [ 0, \infty ) $ \\

When it comes to stability, another of the main {\it troubles} with the Lax-Richtmyer Linear Equivalence Theorem, Rosinger [1-9], is that it is based essentially on the concepts of {\it boundedness} and {\it completeness}. Indeed, in the proof of the implication "convergent $\Longrightarrow$ stable", the Principle of Uniform Boundedness of Linear Operators in Banach Spaces is critically used, and as well known, that principle is not valid in arbitrary normed spaces, but only in Banach spaces, that is, normed spaces which are complete. \\
As argued, however, in Rosinger [1-9], see also section 5 below, stability, as well as convergence of numerical methods, is in fact related with {\it continuity} and {\it compactness}, rather than boundedness and completness. And this fact becomes obvious starting with the next three definitions, and then, with the corresponding results related to the equivalence between convergence and stability. \\

{\bf Definition 2.3.~( Local Stability )} \\

A nonlinear numerical method (2.1) - (2.2) is called {\it locally stable}, if and only if \\

(2.7)~~~ $ \begin{array}{l}
                      \forall~~ \mbox{regular}~ Z\,' \subseteq Z,~~ T \in ( 0, \infty ),~~
                                                                \mbox{compact}~ K \subseteq X ~: \\ \\
                      \exists~~ L\,', \rho\,' > 0 ~: \\ \\
                      \forall~~ \Delta t > 0,~ n \in \mathbb{N},~ u, v \in K \cap X\,'_{n \Delta t} ~: \\ \\
                      ~~~~ n \Delta t \leq T,~ || v - u || \leq \rho\,' ~~~\Longrightarrow~~~
                                    || C^n_{\Delta t} v - C^n_{\Delta t} u || \leq L\,' || v - u ||
            \end{array} $ \\ \\

{\bf Definition 2.4.~ ( Distant Stability )} \\

A nonlinear numerical method (2.1) - (2.2) is called {\it distantly stable}, if and only if \\

(2.8)~~~ $ \begin{array}{l}
                      \forall~~ \mbox{regular}~ Z\,' \subseteq Z,~~ T \in ( 0, \infty ),~~
                                  \mbox{compact}~ K \subseteq X,~~ \rho > 0 ~: \\ \\
                      \exists~~ L\,'' > 0 ~: \\ \\
                      \forall~~ \Delta t > 0,~ n \in \mathbb{N},~ u, v \in K \cap X\,'_{n \Delta t} ~: \\ \\
                      ~~~~ n \Delta t \leq T,~ || v - u || \geq \rho ~~~\Longrightarrow~~~
                                    || C^n_{\Delta t} v - C^n_{\Delta t} u || \leq L\,'' || v - u ||
            \end{array} $ \\ \\

{\bf Remark 2.1.} \\

Two facts are worth noting here : \\

1) There is a significant difference between {\it local stability} and {\it distant stability} as defined above. Indeed, in the former, nearby vectors $u$ and $v$ are involved, since they are subjected to the condition $|| v - u || \leq \rho\,'$ which appears in (2.7). On the other hand, in the latter concept of stability, the vectors $u$ and $v$ involved are supposed to be sufficiently far from one another, as in (2.8) they are required to satisfy the
condition $|| v - u || \geq \rho$. \\

2) The above concept of distant stability is {\it not} the same with a concept of global stability. Indeed, as pointed out above, the vectors $u$ and $v$ involved in (2.8) cannot be arbitrary, since they are restricted by the condition $|| v - u || \geq \rho$.

\hfill $\Box$ \\

Finally \\

{\bf Definition 2.5.~ ( Stability )} \\

A nonlinear numerical method (2.1) - (2.2) is called {\it stable}, if and only if \\

(2.9)~~~ $ \begin{array}{l}
                      \forall~~ \mbox{regular}~ Z\,' \subseteq Z,~~ T \in ( 0, \infty ),~~
                                  \mbox{compact}~ K \subseteq X ~: \\ \\
                      \exists~~ L\ > 0 ~: \\ \\
                      \forall~~ \Delta t > 0,~ n \in \mathbb{N},~ u, v \in K \cap X\,'_{n \Delta t} ~: \\ \\
                      ~~~~ n \Delta t \leq T ~~~\Longrightarrow~~~
                               || C^n_{\Delta t} v - C^n_{\Delta t} u || \leq L || v - u ||
            \end{array} $ \\ \\

The relation between these three concept of stability is given by the easy to prove \\

{\bf Proposition 2.1.~ ( Stability = Local + Distant Stability )} \\

A nonlinear numerical method (2.1) - (2.2) {\it stable}, if and only if it is both {\it locally stable} and {\it distantly stable}. \\

{\bf Remark 2.2.} \\

The need for considering {\it regular} subsets $Z\,' \subseteq Z$ in the above three concepts of stability is due to the fact that even in simple cases, the
set $Z$ in (1.7) fails to be closed, as for instance in the example in Remark 1.1., where \\

$~~~ Z = \{\, ( t, u ) \in [ 0, \infty ) \times \mathbb{R} ~~|~~ t u < 1 \,\} $ \\ \\

{\bf 3. Stability Characterization of the Convergence of \\
        \hspace*{0.45cm} Numerical Methods} \\

{\bf Definition 3.1.} \\

Given a well posed nonlinear semigroup (1.1) - (1.7). A numerical methods (2.1) - (2.2) is called {\it convergent} to that nonlinear semigroup, if and only if \\

(3.1)~~~ $ \begin{array}{l}
                      \forall~~ \mbox{regular}~ Z\,' \subseteq Z,~~ T \in ( 0, \infty ),~~
                                  \mbox{compact}~ K \subseteq X,~~ \epsilon > 0 ~: \\ \\
                      \exists~~ \theta > 0 ~: \\ \\
                      \forall~~ t \in [ 0, T ],~~ \Delta t > 0,~ n \in \mathbb{N},~ u \in K \cap X\,'_t \cap X\,'_{n \Delta t} ~: \\ \\
                      ~~~~ \Delta t,~ | t - n \Delta t | \leq \epsilon ~~~\Longrightarrow~~~
                               || E ( t ) u - C^n_{\Delta t} u || \leq \epsilon
            \end{array} $ \\ \\

{\bf Theorem 3.1.~ ( Convergence $\Longrightarrow$ Distant Stability )} \\

Given a well posed nonlinear semigroup (1.1) - (1.7) and a convergent numerical method (2.1) - (2.2) associated with it. Then that numerical method is distantly
stable. \\

{\bf Proof.} \\

Assume that, on the contrary, there exists a regular subset $Z\,' \subseteq Z,~ T \in ( 0, \infty )$, a compact $K \subseteq X$ and $\rho > 0$, such that
for every $j \in \mathbb{N}$, there exists $\Delta t_j > 0,~ n_j \in \mathbb{N}$, with $n_j \Delta t_j \leq T$, as well as $u_j, v_j \in K \cap
X\,'_{n_j \Delta t_j}$ which satisfy \\

(3.2)~~~ $ || C^{n_j}_{\Delta t_j} v_j - C^{n_j}_{\Delta t_j} u_j || > j || v - u || $ \\

(3.3)~~~ $ || v_j - u_j || \geq \rho $ \\

However, we note that the set \\

(3.4)~~~ $ H = \{\, C^n_{\Delta t} u ~~|~~ \Delta t \in [ 0, T ],~ n \in \mathbb{N},~ n \Delta t \leq T,~ u \in K \cap X\,'_{n \Delta t} \,\} $ \\

is precompact. Indeed, given $\epsilon > 0$, then (3.1) implies that the set \\

(3.5)~~~ $ H_\epsilon = \{\, C^n_{\Delta t} u ~~|~~
                      \Delta t \in [ 0, \theta ],~ n \in \mathbb{N},~ n \Delta t \leq T,~ u \in K \cap X\,'_{n \Delta t} \,\} $ \\

is contained in an $\epsilon$-neighbourhood of \\

(3.6)~~~ $ K\,' =\ \{\, E ( t ) u ~~|~~ t \in [ 0, T ],~ u \in K \cap X\,'_t \,\} $ \\

But \\

(3.7)~~~ $ K\,'$ is compact \\

since in view of (2.6), it is the image of the compact set $Z\,' \cap ( [ 0, T ] \times K )$ under the continuous mapping (1.6). \\

Further, for every given $n \in \mathbb{N}$, the set \\

(3.8)~~~ $ H\,'_n = \{\, C^n_{\Delta t} u ~~|~~ \Delta t \in [ 0, T ],~ u \in K \,\} $ \\

is compact, since it is the image of the compact $[ 0, T ] \times K$ under $n$ iterates of the continuous mapping (2.1) - (2.2). \\

We denote now \\

(3.9)~~~ $ H\,'_\epsilon = \bigcup H\,'_n $ \\

where the union is taken over all $n \in \mathbb{N}$ for which $n \theta \leq T$. And then obviously \\

(3.10)~~~ $ H \subseteq H_\epsilon \cup H\,'_\epsilon $ \\

thus the relations (3.5) - (3.10) complete the proof of (3.4). \\

However, (3.4) implies that $H$ is bounded, therefore, the relations (3.2) and (3.3) contradict one another.

\hfill $\Box$ \\

In a certain sense, a converse of Theorem 3.1. is given below by Theorem 3.2., for which we first need \\

{\bf Definition 3.2.} \\

Given a nonlinear semigroup (1.1) - (1.7). A numerical methods (2.1) - (2.2) is called {\it consistent} with that nonlinear semigroup, if and only if \\

(3.11)~~~ $ \begin{array}{l}
                      \forall~~ \mbox{regular}~ Z\,' \subseteq Z,~~ T \in ( 0, \infty ),~~
                                  \mbox{compact}~ K \subseteq X,~~ \epsilon > 0 ~: \\ \\
                      \exists~~ \delta > 0 ~: \\ \\
                      \forall~~ t \in [ 0, T ],~~ \Delta t > 0,~ u \in K \cap X\,'_{t + \Delta t} ~: \\ \\
                      ~~~~ \Delta t \leq \delta ~~~\Longrightarrow~~~
                               || C_{\Delta t} E ( t ) u - E ( \Delta t ) E ( t ) u || \leq \epsilon \Delta t
            \end{array} $ \\ \\

{\bf Theorem 3.2.~ ( Local Stability $\Longrightarrow$ Convergence )} \\

Given a well posed nonlinear semigroup (1.1) - (1.7) and a locally stable numerical method (2.1) - (2.2) consistent with it. Then the numerical method is
convergent to the nonlinear semigroup. \\

{\bf Proof.} \\

Let be given any regular $Z\,' \subseteq Z,~ T \in ( 0, \infty ),~\mbox{compact}~ K \subseteq X$, and $\epsilon > 0$. Further, let us take any $t \in
[ 0, T ],~ \Delta t > 0,~ n \in \mathbb{N}$, with $n \Delta t \leq T$, and $u \in K \cap X\,'_t \cap X\,'_{n \Delta t}$. \\

We note that \\

(3.12)~~~ $ C^n_{\Delta t} u - E ( t ) u = ( C^n_{\Delta t} u - E ( n \Delta t ) u ) + ( E ( n \Delta t ) u - E ( t ) u ) $ \\

while \\

(3.13)~~~ $ \begin{array}{l}
                      C^n_{\Delta t} u - E ( n \Delta t ) u = \\ \\

                      = \sum_{1 \leq p \leq n} \{ C^{n-p}_{\Delta t} C_{\Delta t} E ( ( p - 1 ) \Delta t ) u -
                                  C^{n-p}_{\Delta t} E ( \Delta t ) E ( ( p - 1 ) \Delta t ) u \}
             \end{array} $ \\

\medskip
However, the set \\

(3.14)~~~$ K\,'' = \{\, C_{\Delta s} E ( s ) v ~~|~~ s, \Delta s \in [ 0, T ],~ v \in K \cap X\,'_s \,\} $ \\

is compact, since it is the image through the continuous mapping (2.1) - (2.2) of $[ 0, T ] \times K\,'$ which in view of (3.7) is compact. Furthermore, in view of (2.2), we have \\

(3.15)~~~ $ K\,' \subseteq K\,'' $ \\

And now, the assumption of local stability applied to $Z\,',~ T$ and $K\,''$ yields $L\,', \rho\,' > 0$, such that for each $1 \leq p \leq n$, we have \\

(3.16)~~~ $ \begin{array}{l}
                       || C^{n-p}_{\Delta t} C_{\Delta t} E ( ( p - 1 ) \Delta t ) u -
                                  C^{n-p}_{\Delta t} E ( \Delta t ) E ( ( p - 1 ) \Delta t ) u || \leq \\ \\
                       ~~~~~~ \leq L\,' \, || C_{\Delta t} E ( ( p - 1 ) \Delta t ) u -
                                                 E ( \Delta t ) E ( ( p - 1 ) \Delta t ) u ||
             \end{array} $ \\

provided that \\

(3.17)~~~ $ C_{\Delta t} E ( ( p - 1 ) \Delta t ) u,~ E ( \Delta t ) u,~
              E ( ( p - 1 ) \Delta t ) u \in K\,'' \cap X\,'_{(n-p) \Delta t} $ \\

and \\

(3.18)~~~ $ || C_{\Delta t} E ( ( p - 1 ) \Delta t ) u -
                                                 E ( \Delta t ) E ( ( p - 1 ) \Delta t ) u || \leq \rho\,' $ \\

Here, the fact that (3.17) is satisfied follows easily from the assumption that $u \in X\,'_{n \Delta t}$, as well as (3.15) and (2.4), (2.5). \\
As for (3.18), we note the following. The assumption of consistency in (3.11) applied to $Z\,',~ T,~ K$ and $\epsilon$ yields $\delta > 0$, such that for $1 \leq p \leq n$, we have \\

(3.19)~~~ $ || C_{\Delta t} E ( ( p - 1 ) \Delta t ) u -
                          E ( \Delta t ) E ( ( p - 1 ) \Delta t ) u || \leq \epsilon \Delta t $ \\

as soon as  $\Delta t \leq \delta$. However, we can obviously assume that \\

(3.20)~~~ $ \Delta t \leq \delta ~~~\Longrightarrow~~~ \epsilon \Delta t \leq \rho\,' $ \\

and then (3.18) holds. \\

Now (3.13), ( 3.16) and (3.19) give \\

(3.21)~~~ $ \Delta t \leq \delta ~~~\Longrightarrow~~~ || C^n_{\Delta t} u - E ( n \Delta t ) u ||
                                       \leq n L\,' \epsilon \Delta t \leq L\,' T \epsilon $ \\

while the assumption of well posedness applied to the compact $Z\,' \cap ( [ 0, T ] \times K )$ results in the uniform
continuity of the corresponding restriction of the mapping (1.6), thus leading to $\delta\,' > 0$, such that \\

(3.22)~~~ $ | t - n \Delta t | \leq \delta\,' ~~~\Longrightarrow~~~
                                     || E ( n \Delta t ) u - E ( u ) || \leq \epsilon $ \\

Finally, (3.21), (3.22) give \\

$~~~ || C^n_{\Delta t} u - E ( t ) u || \leq ( 1 + L T ) \epsilon $ \\

as soon as \\

$~~~ \Delta t \leq \delta,~~~ | t - n \Delta t | \leq \delta\,' $ \\

and the proof of (3.1) is completed.

\hfill $\Box$ \\

{\bf Remark 3.1.} \\

As shown in Rosinger [2, pp. 67-69], Theorems 3.1. and 3.2. above contain as a particular case the Lax-Richtmyer Linear Theorem in its somewhat stronger form as far as the requirements involved in that result are concerned. \\ \\

{\bf 4. Conditions for Nonlinear Equivalence between \\
        \hspace*{0.45cm} Convergence and Stability} \\

The above Theorems 3.1. and 3.2. show that, under the respective conditions, we have the implications \\

~~~~~~"local stability $\Longrightarrow$ convergence $\Longrightarrow$ distant stability" \\

therefore, in view of Proposition 2.1., it follows that \\

~~~~~~"stability $\Longrightarrow$ convergence" \\

In this way, in order to obtain an {\it equivalence} result between convergence and stability, that is, a stability type {\it characterization} of
convergence, it suffices to find conditions - denoted by $({\cal C})$ in the sequel - such that for numerical methods (2.1) - (2.2), the following
implication holds \\

(4.1)~~~ $ \left ( \begin{array}{l}
                                ~~~*)~~~ C_{\Delta t} ~~\mbox{satisfies}~ ({\cal C}) \\ \\
                                ~**)~~~ C_{\Delta t} ~~\mbox{distantly stable}
                    \end{array} ~\right ) ~~~\Longrightarrow~~~ \left (~~~ C_{\Delta t} ~~\mbox{locally stable} ~~~\right ) $ \\

\medskip
And then, one easily obtains \\

{\bf Theorem 4.1.~ ( Nonlinear Equivalence)} \\

Given a well posed nonlinear semigroup (1.1) - (1.7) and a numerical method (2.1) - (2.2) which is consistent with it and also satisfies a condition
of type $({\cal C})$. Then the numerical method is convergent to the nonlinear semigroup, if and only if it is stable.

\hfill $\Box$ \\

{\bf Remark 4.1.} \\

The conditions of type $({\cal C})$ in (4.1) are obviously bridging the {\it gap} between distant stability and local stability. And as seen above, such
a gap occurs regardless of what seems as a rather natural approach to defining the concept of stability of a numerical method associated with a nonlinear
semigroup.

\hfill $\Box$ \\

Here we shall present an {\it explicit} formulation of one specific condition of type $({\cal C})$. Needless to say, there may be many other such
formulations, some of them more useful in applications than other ones. \\

In condition (2.8) defining distant stability, let us, for given $Z\,', T, K$, and $\rho$ denote by \\

(4.2)~~~ $ L\,'' ( Z\,', T, K, \rho ) $ \\

the smallest $L\,''$ for which that relation holds. Then obviously $L\,'' ( Z\,', T, K, \rho )$ is a decreasing function of $\rho$, for $Z\,', T$ and $K$ fixed. Therefore, one obtains easily \\

{\bf Proposition 4.1.} \\

Given a numerical method (2.1) - (2.2) which is distantly stable. Then this numerical method is also locally stable, if and only if \\

(4.3)~~~ $ \begin{array}{l}
                      \forall~~ \mbox{regular}~ Z\,' \subseteq Z,~~ T \in ( 0, \infty ),~~
                                  \mbox{compact}~ K \subseteq X ~: \\ \\
                      ~~~~ \lim_{\,\rho \downarrow 0} L\,'' ( Z\,', T, K, \rho ) < \infty
            \end{array} $

\hfill $\Box$ \\

In view of (4.1), we obtain \\

{\bf Corollary 4.1.~ ( A condition of type $({\cal C})$ )} \\

The above condition (4.3) is a condition of type $({\cal C})$ for numerical methods (2.1) - (2.2). \\ \\

{\bf 5. Continuity and Compactness, instead of Completeness \\
        \hspace*{0.45cm} and Boundedness Are Relevant to Convergence \\
        \hspace*{0.45cm} and Stability} \\

Given a well posed nonlinear semigroup $E ( t )$ in (1.1) - (1.7), and a numerical method $C_{\Delta t}$ in (1.2), (2.2) consistent with it. \\

Let us start by recalling in the {\it particular linear} case, the case in which the Lax-Richtmyer Equivalence Theorem alone operates, the respective {\it linear stability} condition on the linear numerical method, valid on a given time interval $[ 0, T ]$, namely \\

(5.1)~~~ $ \begin{array}{l}
                    \exists~~ L > 0 ~: \\ \\
                    \forall~~ \Delta t > 0,~ n \in \mathbb{N} ~: \\ \\
                    ~~~~ n \Delta t \leq T ~~~\Longrightarrow~~~ || C^n_{\Delta t} || \leq L
            \end{array} $ \\

This simply means that the infinite family of linear operators \\

(5.2)~~~ $ (~ C^n_{\Delta t} ~~|~~ \Delta t > 0,~ n \in \mathbb{N},~ n \Delta t \leq T ~) $ \\

has to be {\it uniformly bounded} on the normed space $( X, ||~|| )$. \\

Second, and as mentioned above, see also Rosinger [1-9], the implication "convergent $\Longrightarrow$ stable" in the Lax-Richtmyer Linear Equivalence Theorem makes essential use of the Principle of Uniform Boundeedness of Linear Operators on Banach Spaces. Thus, as is well known, that theorem {\it cannot} operate in arbitrary normed spaces $( X, ||~|| )$, but only on Banach spaces, that is, normed spaces which are {\it complete}. \\

This is, therefore, the way boundedness and completeness are essentially involved in the Lax-Richtmyer Linear Equivalence Theorem. And as mentioned, see also Rosinger [1-9], that fact leads to major troubles when nonlinear extensions of that theorem are attempted. \\

On the other hand, a simple look at the phenomena involved in the approximation of semigroups by numerical methods does immediately reveal that, instead of boundedness and completeness, two other topological features are in fact essentially involved, namely, {\it continuity} and {\it compactness}. And that fact, as seen above, as well as in Rosinger [1-9], can help in finding {\it nonlinear} generalizations of the Lax-Richtmyer Linear Equivalence
Theorem. \\

Let us now present some of the main related details. \\

As is well known, a vast and particularly important case of approximation of semigroups by numerical methods occurs when the respective semigroups correspond to classical solutions of nonlinear evolution PDEs. And typically in such cases, the respective PDEs are known to have exact solutions, since otherwise, their numerical approximation would of course be pointless. Furthermore, a most important case of such known to exist exact solutions is when, in addition, they are also known to be classical, that is, sufficiently smooth. \\

Consequently, the spaces $X$ in which the respective nonlinear semigroups operate, see (1.1) \\

(5.3)~~~ $ E ( t ) : X_t \longrightarrow X,~~~ t \in [ 0, \infty ) $ \\

can be assumed as given by spaces of sufficiently smooth functions on suitable Euclidean domains, the domains on which the respective nonlinear evolution PDEs are defined. Thus, when endowed with typical norms, such spaces $X$ turn out not to be complete. Furthermore, by completing these spaces $X$ in such norms, one enlarges them considerably, by adding to them large sets of highly nonsmooth functions. \\

On the other hand, having assumed the existence of classical, that is, smooth enough solutions which are elements of such spaces $X$ of smooth enough functions, there is simply no need to further extend those spaces, and in particular, there is no need to complete them in any way. \\

And in fact, even in the general case when the nonlinear semigroups would not necessarily be associated with evolution PDEs, their definition in (1.1) - (1.7) does not in any way require the completeness of the normed spaces $( X, ||~|| )$ on which they are defined. \\

There has also been a rather different line of argument which is claimed to motivate the need for dealing with complete, that is, Banach spaces $( X, ||~|| )$. Namely, it is claimed that the inevitable presence of round-off errors when implementing numerical methods on digital computers leads to that necessity. \\
As it happens, however, such arguments prove to be fallacious, see Rosinger [2, pp. 231-242], [4,5,7,9]. Indeed, the way which is typically suggested for the inclusion of the essentially nonlinear effects of propagation of round-off errors is so unrealistic that it leads to strictly better results in the case of numerical solutions of PDEs, than the well known best possible ones in the case of the numerical solutions of ODEs, see Rosinger [2,4,5,7,9], Isaacson-Keller. \\

So much, therefore, for the requirement of competeness in the general nonlinear situation, that is, for the normed spaces $( X, ||~|| )$ to be in fact Banach spaces, a requirement which gets imposed only by the way the implication
"convergent $\Longrightarrow$ stable" is proved in the particular linear case of the Lax-Richtmyer Equivalence
Theorem. \\

As for {\it continuity}, this condition in one or another form is typically required both on the nonlinear semigroups and the numerical methods aimed to approximate them, as seen in (1.6) and (2.1). \\

And once such a continuity is given, the presence of {\it compactness} is immediate. Indeed, for every given time interval $[ 0, T ]$ and initial value $u \in X_T$, the respective trajectory \\

(5.4)~~~ $ \{\, E ( t ) ~~|~~ t \in [ 0, T ] \,\} $ \\

is always a compact subset of $X$, being the image of the compact interval $[ 0, T ]$ through the continuous mapping
$E ( . ) u$. \\

Therefore, the approximation of the nonlinear semigroup $E ( t )$ by the numerical method $C_{\Delta t}$ is in fact but the approximation of the compact set (5.4). \\

By eliminating completeness, one avoids having to deal with large additional sets of elements in the normed spaces $( X, ||~|| )$, this being a considerable advantage. \\
Second, in infinite dimensional normed spaces, typical when the nonlinear semigroups are associated with evolution nonlinear PDEs, compact sets are bounded, but not the other way round. Furthermore, compact sets are far smaller than bounded sets in general. Consequently, formulating stability conditions involving compact, and not bounded sets, lead to far weaker such conditions. \\

In conclusion, there is a whole variety of advantages which result from formulating equivalence theories between stability and convergence in terms of continuity and compactness, and not in terms of completeness and boundedness, see Rosinger [1-9]. \\

\end{document}